\font\small=cmr5

\def\notdivide{\mid\kern -1.4ex\raise 0.5ex\hbox{\small/}\kern 0.5ex}

\def\lneq{\kern 1.0ex
 \raise 0.1ex \hbox{$\leq$ \kern -2.2ex \lower 0.6ex\hbox{\small
/}\kern 1.5ex}}

\def\dcup {\mathinner{\cup \mkern -8.7mu \rlap{\raise 0.6ex\hbox{.}}\mkern 8.7mu}}

\def\ddbigcup{\mathinner{\bigcup \mkern -15.3mu \rlap{\raise 0.7ex\hbox{.}}\mkern 14.9mu}}
\def\<{\langle}
\def\>{\rangle}
\font\large=cmsy10 scaled \magstep5

\font\small=cmr5

\def\dbigcup{\mathinner{\bigcup \mkern -13.2mu \rlap{\raise 0.6ex\hbox{.}}\mkern 14.9mu}}

\def\C{{\cal C}}

\def\R{{\cal R}}

\def\semidirprod{\times\hskip -0.2em\vrule height 4.6pt depth -0.3pt\hskip 0.2em}

\input xy
\xyoption{all}

 \def\square{\hfill\vbox{\hsize=4pt\hrule\line{\vrule height 2pt depth 2pt\hfill\vrule}\hrule}}


\centerline{CONJUGACY DISTINGUISHED SUBGROUPS}\footnote {}{{\it Mathematics Subject Index} (2010): 20E26,  20E05, 20E18} \footnote {}{{\it Keywords}: residually finite groups, conjugacy separability, conjugacy distinguished subgroup, free-by-finite group, one-relator group, Lyndon group, limit group}

\bigskip
\centerline {Luis Ribes and Pavel A. Zalesskii}\footnote {}{The first author gratefully acknowledges the support of an NSERC research grant and the second author the support of CNP$_q$}

\vskip0.5cm
\hfill {\it To the memory of Oleg V.  Mel'nikov}

\vskip0.5cm

\item{} {\bf Abstract.}   Let $\C$ be  a nonempty
class of finite groups closed under taking subgroups, homomorphic
images and extensions.  A subgroup $H$ of an abstract residually
$\C$  group $R$ is said to be {\it conjugacy $\C$-distinguished}
if whenever $y\in R$, then $y$ has a conjugate in $H$ if and only
if the same holds for the images of $y$ and $H$ in every quotient
group $R/N\in \C$ of $R$. We prove that  in a  group having a normal
free subgroup $\Phi$ such that $R/\Phi$ is in $\C$, every finitely generated subgroup is conjugacy
$\C$-distinguished.  We also prove that finitely generated
subgroups of limit groups, of  Lyndon groups and certain one-relator groups are conjugacy distinguished ($\C$ here is the
class of all finite groups).

\vskip0.5cm
\centerline{\bf  1. Introduction}
\bigskip

In this note we are interested in the following property of a subgroup $H$ of an abstract group $R$: whenever $y\in R$,
then $y$ has a conjugate in $H$ if and only if the same holds for the images of $y$ and $H$ in certain finite   quotient groups $R/N$ of $R$ in a specified class.  We begin by specifying which classes of finite groups will be of interest to us.

Let $\C$ be an {\it extension-closed variety of finite groups}, that is, a nonempty collection of finite groups closed under taking subgroups, homomorphic images and extensions of groups in the collection: if $1\to A\to B \to C \to 1$ is an exact sequence of groups such that $A,C \in \C$, then $B\in \C$. For example $\C$ can be the collection of all finite groups or the collection of all finite solvable groups. If $R$ is a group then its {\it pro-$\C$ completion} $R_{\hat \C}$ is defined to be
$$R_{\hat \C}=\lim\limits_{\displaystyle\longleftarrow\atop {N\in {\cal N}_\C}}R/N,$$
where ${\cal N}_\C$ is the collection of  all the normal subgroups $N$ of $R$ such that $R/N\in \C$. Then $R_{\hat \C}$ is pro-$\C$ group, i.e., a compact, Hausdorff,
totally disconnected topological   group such that $R_{\hat \C}/U\in \C$ whenever $U$ is an open normal subgroup of $R_{\hat \C}$.
The {\it pro-$\C$ topology} of $R$ is defined to be the topology on $R$ that makes it into a topological group so that ${\cal N}_\C$ is a fundamental system of
neighborhoods of the identity element $1$ of $R$. This topology is Hausdorff if and only if the natural homomorphism $R\to R_{\hat \C}$ is injective;
if that is the case one says that  $R$ is a {\it residually $\C$} group, and we think of $R$ as being embedded in $R_{\hat \C}$:  $R\le R_{\hat \C}$. Then
if $X$ is a subset of $R$, we denote its topological closure in $R_{\hat \C}$ by $\bar X$.

  We shall use the standard notation for conjugacy: if $x$ and $r$ are elements of a group $R$, then $x^r= r^{-1}xr$; and we set
$x^R= \{x^r \mid r\in R\}$. If $H$ is a subgroup of $R$, ${\rm N}_R(H)=\{r\in R\mid rH=Hr\}$ denotes as usual its normalizer in $R$, and  ${\rm C}_R(H)=\{r\in R\mid rh=hr, \forall h\in H\}$ its centralizer in $R$.

An abstract  group $R$ is called  {\it conjugacy $\C$-separable}
if for any pair of elements $x,y\in R$, these elements are
conjugate in $R$ if and only if their images in every finite
quotient of $R$ which is in $\C$  are conjugate, or equivalently,
if  $x\not =y^r$ for every $r\in R$, then there exists some
$N\triangleleft R$ with $R/N\in \C$ such that $xN\ne y^rN$ for
every $r\in R$. If $R$ is conjugacy $\C$-separable then   it is
residually $\C$.  One easily checks that a  residually $\C$ group
$R$ is conjugacy $\C$-separable if and only if for any pair  of
elements $x,y\in R$, if $x$ and $y$ are conjugate in  $R_{\hat
\C}$, then they are conjugate in $R$: if $x=y^\gamma$, for some
$\gamma\in R_{\hat \C}$, then there exists some $r\in R$ with $x=
y^r$. If $\C$ is the class of all finite groups, we simply  write
  {\it conjugacy separable}, rather than conjugacy $\C$-separable.

A subgroup $H$ of an abstract residually $\C$  group $R$ is said
to be {\it conjugacy $\C$-distinguished} if whenever $y\in R$,
then $y$ has a conjugate in $H$ if and only if the same holds for
the images of $y$ and $H$ in every   quotient group $R/N\in \C$ of
$R$, or equivalently, $y^R\cap H=\emptyset$ if and only if
$y^{R_{\hat \C}}\cap \bar H=\emptyset$. This means that the
conjugacy class $\bigcup_{r\in R}H^r$ is closed in the pro-$\C$
topology of $R$. If $\C$ is the variety of all finite groups, we
simply write {\it conjugacy distinguished}.

An abstract  group $R$ is said to be {\it free-by-$\C$} if it contains a normal free abstract subgroup $\Phi$ such that $R/\Phi\in \C$.

\medskip
We prove the following results.

\medskip
\noindent {\bf  Theorem A.} {\it
 Let $R$ be a  free-by-$\C$ abstract group  and let $H$ be  a  finitely generated subgroup of $R$ which is closed in its pro-$\C$ topology.
  Then $H$ is conjugacy $\C$-distinguished.}

\bigskip
\noindent {\bf  Theorem B.} {\it Let $R=\langle a_1,...,a_n\mid W^n\rangle$ be a one relator group
with $n>|W|$. Then every finitely generated subgroup $H$ of $R$ is
conjugacy distinguished.}

\medskip
A group $G$ is called fully residually free if for any finite
subset $X$ of $G$ there exists a homomorphism $G$ to a free group
$F$ whose restriction to $X$ is injective. A finitely generated
fully residually free group is called a limit group. Limit groups
have been studied extensively over the last ten years and they
played a crucial role in the solution of the Tarski problem.

\bigskip
\noindent {\bf  Theorem C.} {\it Let $R$ be a limit group and $H$  a finitely
generated subgroup of $R$. Then $H$ is conjugacy distinguished. In
particular, a finitely
generated subgroup of a surface group is conjugacy distinguished.}

\bigskip
The special  case of  Theorem C for  surface groups follows also
from Theorem 1.4 of [3].

Studying equations in free groups Lyndon  [16] introduced groups
$F^{{\bf Z}[t]}$ (later called Lyndon groups) and proved that
these groups are fully residually free; hence a finitely generated
subgroup of a Lyndon group is a limit group.  Lyndon groups play a
very important role in algebraic geometry over groups.
Kharlampovich and Myasnikov  [14] proved  conversely that every
limit group is embeddable into a Lyndon group.

\bigskip
\noindent {\bf  Theorem D.} {\it Let $H$ be a finitely generated subgroup of a Lyndon group $F^{{\bf Z}[t]}$, where $F$ is a free group of arbitrary rank. Then $H$ is
conjugacy distinguished in $F^{{\bf Z}[t]}$.}

\medskip
\noindent{\bf Acknowledgements.} The authors would like to thank
Ashot Minasyan for very useful discussions on Section 3 that led
to considerable improvement of the section.

\vskip0.5cm
\centerline{\bf 2. Free-by-finite groups}
 \bigskip
 In this section we prove that every finitely generated subgroup of a free-by-$\C$
 group which is closed in its  pro-$\C$ topology  is conjugacy $\C$-distinguished.

\bigskip

\noindent {\bf 1. Lemma}  {\it
\medskip
\item {\rm (a)} Let $R$ be a residually $\C$ abstract group
endowed with its pro-$\C$ topology and let   $H$ be a subgroup
of $R$. Assume that    the pro-$\C$ topology of $H$ coincides with
the topology induced from $R$, i.e., that $\bar H= H_{\hat \C}$.
Let $U$ be an open normal subgroup of $R$. Then
$$\overline {U\cap H}= \bar U \cap \bar H.$$
\smallskip
\item {\rm (b)} Let R be a free-by-$\C$   abstract group endowed
with its pro-$\C$ topology. Let $H$ be a finitely generated closed
subgroup and let $U$ be an open normal subgroup of $R$. Then
$$\overline {U\cap H}= \bar U \cap \bar H.$$ }

\smallskip
\noindent {\it Proof.}
\medskip
(a) Note that $\overline{UH}= \bar U\bar H= (UH)_{\hat \C}$ (cf.
[26], Lemma 3.1.4) and  $[R:UH]=
 [R_{\hat \C}:\bar U \bar H]$ (cf. [26], Proposition 3.2.2). So $[UH:U]= [\bar U \bar H:\bar U]$. Therefore
$[H:U\cap H]= [\bar H: \bar U\cap \bar H]$.

Since $\bar H= H_{\hat \C}$, we deduce that if $X\subseteq H$,
then the notation $\bar X$ is unambiguous: it is the closure of
$X$ both in $H_{\hat \C}$ and in  $R_{\hat \C}$. Since $U\cap H$
is open in $H$, we can apply again Proposition 3.2.2 in [26] to
get that $[\bar H: \overline {U\cap H}]= [H: U\cap H]$. Therefore,
$[\bar H: \overline {U\cap H}]=[\bar H: \bar U\cap \bar H]$. Since
$\overline {U\cap H}\le \bar U\cap \bar H$, we deduce that
$\overline {U\cap H}=  \bar U\cap \bar H$.

\smallskip

(b) We shall prove that in this case $\bar H= H_{\hat \C}$, and
then the result will follow from part (a). Let  $\Phi$ be a free
normal subgroup of $R$ such that $R/\Phi\in \C$. The pro-$\C$
topology of $\Phi$ coincides with the topology induced from the
pro-$\C$ topology of $R$ (cf. [26], Lemma 3.1.4 (a)). Since
$H/(\Phi \cap H)\in \C$, the subgroup $\Phi \cap H$ is open  in
the  pro-$\C$ topology of $H$ as well as in the topology induced
from $\Phi$. Hence it suffices to show that the    pro-$\C$
topology of $\Phi \cap H$  coincides with the topology induced
from $\Phi$. Since $\Phi \cap H$ is a closed finitely generated
subgroup of the free group $\Phi$,  there exists an open subgroup
$W$ of $\Phi$ such that $W= L* H$, for some subgroup $L$ of $W$
(cf. [24], Lemma 3.2).  Therefore the pro-$\C$ topology of $H$
coincides with that induced from the pro-$\C$ topology of $W$ (cf.
[26],   Corollary 3.1.6 (a)), and hence with the topology  induced
from the pro-$\C$ topology of $R$, since $W$ is open in $R$, i.e.,
$\bar H= H_{\hat \C}$ as needed.
  \hfill~$\square$

\bigskip
The  last statement of the lemma above is a special case of a more
general result ([27],  Proposition 2.3), where one does not
require that $U$ be open.


 The next two   results sharpen Lemma 2.2 and  Theorem 3.2 in [27], where they are  proved only for finitely generated groups $R$.

\bigskip

\noindent {\bf   2. Lemma} {\it    Let  $H\in \C$ be a group of
prime order $p$. Let $R=\Phi\semidirprod H$ be a semidirect
product, where $\Phi$ is an abstract free group.  Then  there is a
free factor $\Phi_1$  of $\Phi$   such that

\medskip
\item{\rm (a)} ${\rm N}_R(H)  =H\times \Phi_1$   and ${\rm N}_{R_{\hat
\C}}(\bar H)=  H\times  (\Phi_1)_{\hat \C}$;

\smallskip
\noindent and

\smallskip
\item{\rm (b)} ${\rm C}_\Phi(H) = \Phi_1$  and ${\rm C}_{\Phi_{\hat
\C}}(H) = (\Phi_1)_{\hat \C} $.

\medskip
\noindent Consequently,
\medskip

\item{\rm (a$'$)} $\overline {{\rm N}_R(H)}=
{\rm N}_{R_{\hat \C}}(H)$;

\smallskip
\item {\rm (b$'$)} $\overline {{\rm C}_\Phi(H)}=
{\rm C}_{\bar \Phi}(H)$.}

\medskip

\noindent {\it Proof.}  By a theorem of   Dyer-Scott (cf. [6],
Theorem 1)  the group $R$ is a free product
$$R=\big[
\hbox{\textfont2=\large$*$} _{i\in I} (C_i\times \Phi_i)\big]*L,$$
 where  $L$ and each $\Phi_i$ are free groups   and the $C_i$ are groups of
order $p$. Since every finite subgroup of $R$ of order $p$ is
conjugate to one of the $C_i$  (cf.   [18], Corollary 4.1.4), we
may assume without loss of generality that $H=C_{i_1}$, for some
fixed $i_1\in I$. Then $R=  (C_{i_1}\times \Phi_{i_1}) * R_1=
(H\times \Phi_1) *R_1$, where $\Phi_1= \Phi_{i_1}$ and
$$R_1=    \big[ \hbox{\textfont2=\large$*$} _{i\in I-\{i_1\}} (C_i\times \Phi_i)\big]*L.   $$
It follows that ${\rm N}_R(H) =H\times \Phi_1$ (cf. [18],
Corollary 4.1.5), and since $H$ is abelian, ${\rm C}_R(H)=H\times
\Phi_1$. Hence $\Phi_1=  {\rm N}_\Phi(H)= {\rm C}_\Phi(H)\le
\Phi$. Now, ${\rm C}_\Phi(H)$  is the subgroup of fixed points of
$\Phi$ under the action of $H$, and so ${\rm C}_\Phi(H) = \Phi_1$
is a free factor of $\Phi$ (cf. [6], Theorem 2).   This implies
that $( \Phi_1)_{\hat \C}= \overline {\Phi_1}$ (cf. [26],
Corollary 3.1.6).

Finally observe that $R_{\hat \C}= (H\times (\Phi_1)_{\hat
\C})\amalg (R_1)_{\hat \C}$ (the free pro-$\C$ product);  so ${\rm
N}_{R_{\hat \C}}(H)=  H\times (\Phi_1)_{\hat \C}$ (cf. [26],
Theorem 9.1.12), and since $H$ is abelian,  ${\rm C}_{R_{\hat
\C}}(H)=  H\times (\Phi_1)_{\hat \C}$. Therefore,  ${\rm
N}_{R_{\hat \C}}(H)=  H\times \overline {\Phi_1}= \overline {{\rm
N}_R(H)}$ and  ${\rm C}_{R_{\hat \C}}(H)=  H\times \overline
{\Phi_1}= \overline {{\rm C}_R(H)}$. Thus, ${\rm C}_{\bar
\Phi}(H)={\rm C}_{R_{\hat \C}}(H)\cap \overline\Phi= (H\times
\overline {\Phi_1})\cap \overline\Phi=\overline {\Phi_1}=
\overline {{\rm C}_\Phi(H)} $.
 This concludes the proof of all parts of the lemma.
   \hfill$\square$

\bigskip
\noindent {\bf 3. Theorem}  {\it A
free-by-$\C$ group  $R$ is conjugacy $\C$-separable.}

\medskip \noindent{\it Proof.}  To fix the
notation, say  that $\Phi\triangleleft R$, where $\Phi$  is an
abstract free group such that $R/\Phi \in \C$.  Then we may think of  $\Phi_{\hat \C}$ as   an open subgroup of
 $R_{\hat \C}$ (cf. [26], Lemma 3.1.4 (a)).
Let $x,y\in R$ and let $x^\gamma=y$, where $\gamma\in R_{\hat
\C}$.  We have to show that $x$ and $y$ are conjugate in $R$. We
may assume that $x\not=1$. Since $R_{\hat \C}=R\Phi_{\hat \C}$,
we have $\gamma=r\eta$,  for some $r\in R$, $\eta\in \Phi_{\hat
\C}$.  So replacing $x$ by $x^r$ and $\gamma$ by $\eta$,  we may
assume that $\gamma$ is in $\Phi_{\hat \C}$. Then $y\in \langle
x\rangle \Phi_{\hat \C}\cap R= \langle x\rangle \bar  \Phi \cap R=
\langle x\rangle (\bar  \Phi \cap R)=\langle x\rangle \Phi$.
Hence from now on we may also assume that  $R= \langle x\rangle
\Phi$. Note that $R_{\hat \C}= \langle x\rangle \Phi_{\hat \C}$.
Since  $R_{\hat \C}/\Phi_{\hat \C}$ is abelian, we have
$x^{-1}\gamma^{-1} x\gamma\in \Phi_{\hat \C}$, i.e., $x\Phi_{\hat
\C}= x^\gamma\Phi_{\hat \C}$.  On the other hand we have that the
natural map $\rho: R/\Phi\longrightarrow R_{\hat \C}/\Phi_{\hat
\C}$ is a bijection.  Since $\rho(y\Phi)= y\Phi_{\hat \C}=
x^\gamma\Phi_{\hat \C}=  x\Phi_{\hat \C}= \rho(x\Phi)$, we deduce
that $y\Phi= x\Phi$. So from now on we assume that
$$R=  \langle x\rangle \Phi, \quad y=x^\gamma\in R,\ {\rm with}\  \gamma\in \Phi_{\hat \C}, \ {\rm and}\  y\Phi= x\Phi.\eqno{(1)}$$
 Now we distinguish two cases.
\medskip
\noindent {\it Case 1.} The order of  $x$ is infinite. Let $n$ be
a positive integer such that $x^n\in \Phi $.  So $y^n\in \Phi$ and
$y^n= (x^n)^\gamma$. According to a result of Baumslag and Taylor
(cf. [17], Proposition 4.8), free groups are conjugacy
$\C$-separable. We deduce that    $y^n$  and $x^n$ are conjugate
in $\Phi$. Say  $f^{-1}x^n f= y^n$, where $f\in \Phi$. Replacing
$x$ with $fxf^{-1}$,  we may assume that $y^n= x^n$. Therefore
$\gamma\in {\rm C}_{\Phi_{\hat \C}}(x^n)= {\rm C}_{R_{\hat
\C}}(x^n)\cap \Phi_{\hat \C}$.  Write  $\Phi= \Phi_1*\Phi_2$ ,
where $\Phi_1$ is a free subgroup of $\Phi$ of finite rank such
that $x^n\in \Phi_1$. Then  $\Phi_{\hat \C}= ( \Phi_1)_{\hat
\C}\amalg (\Phi_2)_{\hat \C}= \overline { \Phi_1}\amalg \overline
{ \Phi_1}$, the free pro-$\C$ product (here we are using Corollary
3.1.6 in [26]). Note that ${\rm C}_\Phi (x^n)= {\rm C}_{\Phi
_1}(x^n)$ and ${\rm C}_{\Phi_{\hat \C}}(x^n)= {\rm
C}_{(\Phi_1)_{\hat \C}}(x^n)$ (cf. [26], Theorem 9.1.12). Since
$\Phi_1$ has finite rank we can use Corollary 2.8 in [27] to get
that $\overline{{\rm C}_{\Phi_1}(x^n)}= {\rm C}_{\overline
{\Phi_1}}(x^n)$, and so $\overline{{\rm C}_{\Phi}(x^n)}= {\rm
C}_{\overline {\Phi}}(x^n)=  {\rm C}_{ \Phi_{\hat \C}}(x^n)$.
Therefore, $\gamma\in   \overline{{\rm C}_{\Phi}(x^n)}$.  Since
${\rm C}_{\Phi}(x^n)\le  {\rm C}_{R}(x^n)$, we have
$\overline{{\rm C}_{\Phi}(x^n)}\le \overline{{\rm C}_{R}(x^n)}$.
Hence  $\gamma\in   \overline{{\rm C}_{R}(x^n)}$. Thus $x,y,
\gamma\in  \overline{{\rm C}_{R}(x^n)}$.

 Since $x^n\not=1$  and $\Phi$ is free, ${\rm C}_\Phi(x^n)$ is cyclic, say  ${\rm C}_\Phi(x^n)= \langle z\rangle$  and $z^m= x^n$, for some natural number $m$. Using the uniqueness of $m$-th roots in $\Phi$, we get that ${\rm C}_R(x^n)=  {\rm C}_R(z)$. Hence $x\in {\rm C}_R(z)$, i.e., $x$ and $z$ commute.

Since $R= \langle x\rangle\Phi$, we obtain that
${\rm C}_R(x^n)=  \langle x\rangle  {\rm C}_\Phi(x^n)=\langle x\rangle \langle z\rangle$; therefore  ${\rm C}_R(x^n)$  is abelian, and  hence so is  $\overline {{\rm C}_R(x^n)}$.   This implies that $x=y$; thus the result holds in this case.

\medskip
\noindent {\it Case} 2. The order of  $x$ is finite.  Observe that
$\langle x \rangle$ is isomorphic to a subgroup of $R/\Phi$, and
so $\langle x \rangle\in \C$.  We proceed by induction on the
order of  $x$.
\smallskip

 {\it Subcase} 2 (a).  The order of  $x$ is  $p$,  a prime.  Then $y$ is also of order $p$. By a theorem of   Dyer-Scott (cf. [6], Theorem 1)
 the group
$R$ is a free product

$$R= \langle x\rangle \Phi=\big[
\hbox{\textfont2=\large$*$} _{i\in I} (C_i\times \Phi_i)\big]*L,$$
 where  $L$ and each $\Phi_i$ are free groups   and the $C_i$ are groups of
order $p$. Suppose $x$ and $y$ are not conjugate in $R$. Since
every finite subgroup of $R$ of order $p$ is conjugate in $R$ to
one of the $C_i$ (cf. [18], Corollary 4.1.4), we may assume that
$C_{i_1}= \langle x\rangle$ and $C_{i_2}=\langle y\rangle$, where
$i_1, i_2\in I$ and $i_1\ne i_2$. Hence $R= (C_{i_1}\times
\Phi_{i_1})* (C_{i_2}\times \Phi_{i_2})* R_1$, where
$$R_1= \big[
\hbox{\textfont2=\large$*$} _{i\in I-\{i_1, i_2\}}(C_i\times
\Phi_i)\big]*L.$$ Define  $\tilde  R=C_{i_1} * C_{i_2}$ and let
$\varphi: R \to \tilde R$ be a natural epimorphism that sends
$C_{i_1} $ and $C_{i_2}$ identically to their corresponding copies
in $\tilde  R$ and sends $\Phi_{i_1}$, $\Phi_{i_2}$ and $R_1$ to
$1$. Then $x$ and $y$ are not conjugate in the free pro-$\C$
product $\tilde  R_{\hat \C}=  C_{i_1}\amalg C_{i_2}$ (cf. [26],
Theorem 9.1.12). However, the epimorphism $\varphi$ induces an
epimorphism $\hat \varphi: R_{\hat \C}\to \tilde  R_{\hat \C}=
C_{i_1}\amalg C_{i_2}$  (cf. [26], Proposition 3.2.1), and so
$x^{\varphi (\gamma)}=y$ in $\tilde  R_{\hat \C}= C_{i_1}\amalg
C_{i_2}$, a contradiction.

\smallskip

 {\it Subcase} 2 (b).  The order of $x$ is finite but not a prime. Choose a natural number $n$ such that the order of $x^n$ is a prime. By the  subcase 2 (a) above, replacing $x$ by a certain conjugate in $R$, we may assume that   $x^n= y^n$, and so $\gamma$ centralizes $x^n$; hence  $\gamma  \in {\rm C}_{ \Phi_{\hat \C}}(x^n)=  \overline { {\rm C}_{\Phi}(x^n)}$  (the last equality is the content of Lemma 2  (b$'$) above).    Put
$H=\langle x\rangle  {\rm C}_{\Phi}(x^n)$.    Since $x$ normalizes
${\rm C}_{\Phi}(x^n)$,  $H$ is a subgroup of $R$.  By Lemma 2,
${\rm C}_{\Phi}(x^n)$  is a free factor of $\Phi$, and so it is
closed in  $\Phi$. Hence ${\rm C}_{\Phi}(x^n)$ is closed in $R$;
moreover the pro-$\C$ topology on it induced from $R$ is its full
pro-$\C$ topology (cf. Corollary 3.1.6 in [26]).  Since $\langle
x\rangle$ is finite, $H$ is closed in  $R$ and $H_{\hat \C}=  \bar
H$  (this follows from Corollary 3.3 in [24]). Therefore, $H_{\hat
\C}= \bar H= \langle x\rangle \overline { {\rm C}_{\Phi}(x^n)}$.
It follows that $x, y  \in H$  and $\gamma\in H_{\hat \C}$.  Hence
we may assume that $R=H=\langle x\rangle {\rm C}_{\Phi}(x^n)$.
Moreover conditions (1) still hold, where now ${\rm
C}_{\Phi}(x^n)$ plays the role of $\Phi$. Note that then $\langle
x^n\rangle$ is a central subgroup of  $R$, and $R/\langle
x^n\rangle=  (\langle x\rangle/\langle x^n\rangle) {\rm
C}_{\Phi}(x^n)$, where, with a certain abuse of notation,  we
identify ${\rm C}_{\Phi}(x^n)$ with its isomorphic image in
$R/\langle x^n\rangle$.  Denote by  $\tilde x$  and $\tilde y$ the
images of $x$ and $y$ in $R/\langle x^n\rangle$, respectively. So
$R/\langle x^n\rangle= \langle \tilde x\rangle {\rm
C}_{\Phi}(x^n)$. Note that the order of $\tilde  x$ is strictly
smaller than the order of $x$; $\tilde  y= \tilde x^\gamma$,  with
$\gamma \in
 {\rm C}_{\Phi}(x^n)$,  and $(R/\langle x^n\rangle)/ {\rm C}_{\Phi}(x^n)\cong  \langle \tilde  x\rangle\in \C$. By the induction hypothesis, there exists
 some $f\in {\rm C}_{\Phi}(x^n)$  such that $\tilde  y=  \tilde x^f$.  Replacing $x$ with $x^f$  and $\gamma$ with $f^{-1}\gamma$, we may assume that
 $\tilde  y= \tilde x$;   observe that conditions (1) still hold,  with  ${\rm C}_{\Phi}(\langle x^n\rangle)$  playing the role of $\Phi$.   Therefore
  $y= xc$, for some $c\in \langle x^n\rangle$. Since $x{\rm C}_{\Phi}(\langle x^n\rangle)= y{\rm C}_{\Phi}(\langle x^n\rangle)$,  and ${\rm C}_{\Phi}(x^n)$ is a free group, we have $c=1$. Thus $x=y$, and the result follows.
 \hfill $\square$

\bigskip

\noindent {\it  Proof of Theorem A.}    This is equivalent to
proving that if $a\in R$ and $a^\gamma=\gamma^{-1}a\gamma\in \bar
H$, where $\gamma\in R_{\hat \C}$, then there exist $c\in R$ such
that $c^{-1}ac\in H$.

It follows from a result of Scott ([29]) that  $R$ is the
fundamental group $\Pi^{abs}= \Pi^{abs}(\R, \Delta)$ of a graph of
groups $(\R, \Delta)$ over a graph $\Delta$ such that each vertex
group $\R (v)$ is in $\C$ ($v\in V(\Delta)$).

Since  $R$ is  free-by-$\C$, there exists a subgroup   $\Phi$ of
$R$ which is free  and open in the pro-$\C$ topology of $R$.
 The pro-$\C$ topology of $R$ induces on $\Phi$ its own full pro-$\C$ topology ([26], Lemma 3.1.4 (a)).

\medskip
\noindent {\it Case}  (i).  The element $a$ has finite order.
\smallskip
 \noindent The pro-$\C$ topology of $R$ induces on  $H$ its own full pro-$\C$ topology, so that one can make the identification
 $\bar H=H_{\hat \C}$ (indeed since $H$ is finitely generated,  one can write $\Phi=  \Phi_1 *  \Phi_2$, where $\Phi_1$ is a
 free group of finite rank
such that   $H\cap \Phi\le \Phi_1 $;  since $\Phi_1$ is closed in
the pro-$\C$ topology of $R$ and $\overline {\Phi_1}=
(\Phi_1)_{\hat \C}$ by [26], Corollary 3.1.6, so by Lemma 3.1 in
[27],
 $\bar H=H_{\hat \C}$). Observe that $H$ is also a
finitely generated free-by-$\C$ group, and so, using a result of
Karrass, Petrowski and Solitar ([13], Theorem 1) or the result of
Scott mentioned above, $H$ is the fundamental group $\Pi^{'abs}=
\Pi^{abs}(\R', \Delta')$ of a graph of groups in $\C$, $(\R',
\Delta')$, over a finite
 graph $\Delta'$; and $\bar H= H_{\hat \C}$
  is the pro-$\C$ fundamental group of $(\R', \Delta')$.
 In addition  we may make the identification $\R'(v)= \Pi^{'abs}(v)$, a subgroup of $H$, for every vertex
 $v$ of $\Delta'$ ([27], Section 0).
 Since $\gamma^{-1}a\gamma \in \bar H$ has finite order,
it  is conjugate in $\bar H= H_{\hat \C}$ to an element of some
vertex group $\R'(w)= \Pi^{'abs}(w)\le H$ ([32], Theorem 3.10).
Therefore, since $H_{\hat \C}\le R_{\hat \C}$, $a$ is conjugate in
$R_{\hat \C}$ to an element, say $b$, of $H$. Thus, by Theorem  3,
there exists $c\in R$ with $c^{-1}ac= b\in H$.

\smallskip

\noindent {\it Case}  (ii).  The element $a$ has infinite order.
\smallskip
 \noindent   Since $\Phi \cap H$ is  finitely generated
 and it is closed in the pro-$\C$ topology of $\Phi$, it follows from [24],   Lemma 3.2 that there exists an open subgroup $U$ of $\Phi$ containing
 $\Phi\cap H$ such
 that $\Phi=(\Phi\cap H)*L$, for some closed subgroup $L$ of $\Phi$. Replacing $\Phi$ with $U$ we may assume that
$\Phi=(\Phi\cap H)*L$.

Since  $R$ is dense in  $R_{\hat \C}$ and $\bar \Phi$  is open in
$R_{\hat \C}$, we have that $R_{\hat \C}= R\bar \Phi$. So $\gamma=
r\gamma_1$, for some $r\in R, \gamma_1\in \bar \Phi$. Therefore,
replacing $a$ with $r^{-1}ar$, we may assume that $\gamma=\gamma_1
\in \bar \Phi= \Phi_{\hat \C}$. Since $\Phi$ has finite index in
$R$,  we have $1\ne a^n\in \Phi$, for some natural number $n$.
Observe that the pro-$\C$ topology of $R$ induces on $\Phi\cap H$
(respectively, on $L$) its full pro-$\C$ topology ([26], Corollary
3.1.6); therefore,
$$\bar \Phi= \Phi_{\hat \C}=(\Phi\cap H)_{\hat \C}\amalg L_{\hat \C}=\overline {(\Phi\cap H)}\amalg \bar L,$$
 the free pro-$\C$ product ([26], Section 9.1). By Lemma 1, $\overline {\Phi\cap H}= \bar \Phi \cap \bar H$; so  $\gamma^{-1}a^n\gamma\in \overline {(\Phi\cap H)}$.

 We deduce from ([25], Proposition 2.9) that $a^n$ is nonhyperbolic as an element of the free product $\Phi=(\Phi\cap H)*L$; i.e., $a^n$ is conjugate in $\Phi$ to an element of  either $\Phi\cap H$ or $L$; in fact it must be conjugate in $\Phi$ to an element of $\Phi\cap H$,  since otherwise $\overline {\Phi\cap H}$ would contain a conjugate in  $\bar \Phi$ of a nontrivial element of $\bar L$, which is not possible ([26], Theorem 9.1.12).  Say   $c^{-1}a^nc\in \Phi \cap H$, for some $c\in \Phi$.  Then $(\gamma^{-1}c)c^{-1}a^nc(c^{-1}\gamma)\in \overline {\Phi\cap H}$; therefore using again Theorem 9.1.12  in [26], we have that  $c^{-1}\gamma\in \overline {\Phi\cap H}$. We deduce that $c^{-1}\gamma\in \bar H$. Since $\gamma^{-1}a\gamma\in \bar H$, we have
$(\gamma^{-1}c)c^{-1}a c(c^{-1}\gamma)\in \bar H $, and therefore
 $c^{-1}a c\in \bar H$. Now, since $H$  is closed in the pro-$\C$ topology of $R$ by assumption, we have $\bar H\cap R= H$. Thus
$$c^{-1}a c\in \bar H \cap R=H,$$
as needed.\hfill~$\square$

\bigskip

In the profinite topology of a  free-by-finite group every
finitely generated subgroup is closed (this follows easily from
[10],   Theorem 5.1). Hence one deduces the following result
(essentially proved in [5] when the subgroup $H$ is cyclic).
\bigskip
\noindent {\bf 4. Corollary}  {\it  Let $R$ be a  free-by-finite abstract group, and let $H$ be a
finitely generated   subgroup of $R$.  Then $H$ is conjugacy
distinguished}.

\bigskip

\noindent{\bf 5. Remark}
 The condition in Theorem A that $H$ is closed in the pro-$\C$
topology of $R$ is necessary. For example, let $R= {\bf Z}$, the
free group of rank $1$. Let $p$ be a prime number and let $\C$
consist of all finite $p$-groups, so that the pro-$\C$ topology
is, in this case, the pro-$p$ topology. Consider the subgroup $H=
q{\bf Z}$ of ${\bf Z}$, where $q$ is a prime, $q\ne p$.  Then $H$
is not closed in the pro-$p$ topology of ${\bf Z}$, but if
$\varphi_n: {\bf Z}\to {\bf Z}/p^n{\bf Z}$ is the natural
epimorphism ($n=1,2,\dots$), then   $\varphi_n (H)= {\bf
Z}/p^n{\bf Z}$. Therefore, if $x\in {\bf Z}- H$, one has $\varphi
_n(x)\in \varphi_n(H)$ for each $n$. Thus, since  ${\bf Z}$ is
abelian, $H$ is not conjugacy $\C$-distinguished  in  ${\bf Z}$.

\bigskip
\centerline{  \bf 3. Virtual retracts}

\bigskip
In this section we prove Theorems B and C. We first introduce the concept of `virtual retract' for a subgroup $H$ of a group $R$ and show
how this property  helps to establish that $H$  is  conjugacy
distinguished.

\medskip

We say that a subgroup $H$ of an abstract group $R$ is a {\it $\C$-virtual retract} of $R$ if there is an open  subgroup $U$ in the
 pro-$\C$ topology of $R$ such that $H\le U$ and there exists an epimorphism $f: U\to H$ that is the identity when restricted to $H$ (such an epimorphism is called a retraction and the subgroup $H$ is called a retract of $U$); in other words
$U= K\semidirprod H$, where $K$ is a subgroup of $U$. When  $\C$
is the class of all finite groups, we simply say that $H$ is a
{\it virtual retract} of $R$. Observe that if $H$ is a
$\C$-virtual retract of $R$, then $H$ is closed in the pro-$\C$
topology of $R$ and $\bar H= H_{\hat \C}$ ([26], Lemmas 3.1.4 and
3.1.5).

\bigskip
\noindent {\bf 6. Lemma} {\it  Let $R$ be an abstract   group.

\medskip
\item {\rm (a)} Assume that  $R$ is conjugacy $\C$-separable and
let $H$ be a retract of $R$.  Then $H$ is conjugacy \break
$\C$-distinguished.

\smallskip
  \item {\rm (b)} Assume that every open (in the pro-$\C$ topology)
subgroup of $R$ is conjugacy $\C$-separable. Let $H$ be a
$\C$-virtual retract of $R$.  Then $H$ is conjugacy
$\C$-distinguished.
 \item {}}

\medskip

\noindent {\it Proof.}  We first reduce (b) to (a). Let $V$ be an
open subgroup of $R$ containing $H$. Observe that the conjugacy
class of $H$ in $R$ is a finite union of translates of the
conjugacy class  of $H$ in $V$: $\bigcup_{r\in R}
H^r=\bigcup_{i=1}^n(\bigcup_{u\in V}  H^ u)^{t_i}$, where $\{t_1,
\dots, t_n\}$ is a transversal of $V$ in $R$. Assume that the
hypotheses of (b) hold, and let $U$ be an open subgroup of $R$
containing $H$ such that $H$ is a retract of $U$. Recall that $H$
is a $\C$-distinguished subgroup of $R$ if and only if its
$R$-conjugacy class $H^R= \bigcup_{r\in R} H^r$ is closed in the
pro-$\C$ topology of $R$. It follows then from the above
observation that to show that $H$ is $\C$-distinguished as a
subgroup of $R$, it suffices to prove that it is
$\C$-distinguished as a subgroup of $U$. Therefore (b) follows
from (a).

 To prove (a) note first that as pointed out
above  $\bar H= H_{\hat \C}$. Let  $f:R\to H$ be  an epimorphism
that is the identity map on $H$. Consider the extension $\hat
f:R_{\hat \C}\to H_{\hat \C}=\bar H$ of $f$ to the pro-$\C$
completions of $R$ and $H$. Then the restriction of $\hat f$ to
$\bar H$ is the identity map on $\bar H$.  Let $r\in R$,
$\gamma\in R_{\hat \C}$ be such that $r^\gamma\in \bar H$. We need
to show that there exists some $s\in R$ with $r^s\in H$.

One has
$$r^\gamma= \hat f(r^\gamma)= f(r)^{\hat f(\gamma)}.$$
Therefore,
 $r$ and $f(r)$ are conjugate in
$R_{\hat \C}$. Since $R$ is conjugacy $\C$-separable by
assumption, $r$ and $f(r)$ are conjugate in $R$; say $s^{-1}rs =
f(r)$, for some  $s\in R$. So    $s^{-1}rs \in H$, as
needed.\hfill~\square

\bigskip
Next we prove a slightly different version of the lemma above
under   less restrictive assumptions.

\bigskip
\noindent {\bf 7. Proposition} {\it  Let $R$ be an abstract group
endowed with its pro-$\C$ topology and let $H$ be a closed
subgroup of $R$. Assume  that every open   subgroup of $R$ is
conjugacy $\C$-separable. Let $U$ be an open subgroup of $R$ such
that $H\cap U$ is  a  retract of $U$. Suppose ${\rm C}_{
R_{\hat\C}}(h)$ is abelian for every $1\neq h\in H$. Then for
every element $g$ of $R$ of infinite order the equality $g^R\cap
H=\emptyset$ implies $g^{R_{\hat \C}}\cap \overline H=\emptyset$.
In particular, if $R$ is torsion-free, then $H$ is conjugacy
$\C$-distinguished.}

\medskip

\noindent {\it Proof.} As pointed out above, the pro-$\C$ topology
of $U\cap H$ coincides with the topology induced from the pro-$\C$
topology of $U$ (and hence of $R$). Since $U\cap H$ is open in
both topologies of $H$ (its own and the one induced from the
pro-$\C$ topology of $R$), we deduce that the same holds for $H$,
i.e., $\bar H= H_{\hat \C}$.

 Let $g\in R$, $\gamma\in R_{\hat \C}$ be such that
$g^\gamma\in \bar H$. We need to show that there exists some $c\in
R$ with $g^c\in H$. Let   $f:U\to U\cap H$ be an epimorphism which
is the identity map on $U\cap H$. Consider  the extension $\hat
f:\bar U= U_{\hat \C}\to \overline{U\cap H}= (U\cap H)_{\hat \C}$
of $f$ to the pro-$\C$ completions of $U$ and $U\cap H$. Then the
restriction of $\hat f$ to $\overline{U\cap H}$ is the identity
map on $\overline{U\cap H}$.

Since $\bar U$ is open in $R_{\hat \C}$, one has $R_{\hat \C}=
R\bar U$. Write $\gamma= r\gamma_1$, with $r\in R$ and
$\gamma_1\in \bar U$. Then $g^\gamma= (g^r)^{\gamma_1}$. Therefore
replacing $\gamma$ with $\gamma_1$, if necessary, we may assume
that $\gamma\in \bar U$. Let $n$ be a natural number such that
$u=g^n\in U$. Since $u^\gamma= (g^\gamma)^n\in \bar H\cap \bar
U=\overline{H\cap U}$ (see Lemma 1), one has
$$u^\gamma= \hat f (u^\gamma)= f(u)^{\hat f(\gamma)}.\eqno{(2)}$$
Therefore,
 $u$ and $f(u)$ are conjugate in
$\bar U$. Since $U$ is conjugacy $\C$-separable by assumption, $u$
and $f(u)$ are conjugate in $U$; say $(s^{-1}gs)^n = f(g^n)$,
where $s\in U$. So    $(s^{-1}gs)^n \in H$.     Then (2) can be
rewritten as
$$\big((g^s)^n\big)^{s^{-1}\gamma}=   \big(f(g^s)^{n}\big)^{\hat f(s^{-1}\gamma)}$$
since $f(s)$ and $\hat f(s^{-1}\gamma)$ are defined. Replacing $g$
with $g^s=s^{-1}gs$ and $\gamma$ with $s^{-1}\gamma$, we may
assume that $g^n\in H$. Hence the equation above reads as
$(g^n)^\gamma=(g^n)^{\hat f(\gamma)}$. It follows that
$\gamma(\hat f(\gamma))^{-1}\in {\rm C}_{ R_{\hat\C}}(g^n)$, which
is an abelian group by hypothesis, because $g^n\ne 1$. Since we
also have that $g\in {\rm C}_{R}(g^n)$, one gets
$$g^\gamma=  \big(g^{\gamma \hat f(\gamma)^{-1}}\big)^{\hat f(\gamma)}= g^{\hat f(\gamma)}.$$
From  $g^\gamma\in \bar H$, we deduce that $g^{\hat f(\gamma)}\in
\bar H$. Therefore $g\in \bar H$, because $\hat f(\gamma)\in \bar
H$. Since $H$ is closed in the pro-$\C$ topology of $R$, we obtain
that  $g\in \bar H\cap R= H$, so that in this case we can take
$c=1$, proving the result.\hfill~\square

\bigskip

\noindent  {\bf 8. Remarks}
\medskip
8.1  In contrast to residual finiteness and subgroup separability,
the conjugacy separability property is not inherited by subgroups
of finite index. So the assumption  in Lemma 6  that every open
subgroup of the group $R$ is conjugacy $\C$-separable  is
essential.
\smallskip
8.2 The assumption in Proposition 7 that ${\rm C}_{\hat R}(h)$  is
abelian for every $h\in H$  is, in principle, not easily
verifiable. With this in mind we make use of the following result
of Minasyan ([19], Proposition 3.2): a group  $R$ and all its
subgroups of finite index are conjugacy separable if and only if
$R$ is conjugacy separable and ${\rm C}_{\hat R}(g)= \overline
{{\rm C}_R(g)}$, for every $g\in R$  (this result  has been
extended by Ferov in [7], Theorem 4.2, to a corresponding
equivalence for the property of conjugacy $\C$-separability). The
idea is that in certain cases it suffices to know the abelianness
of ${\rm C}_R(g)$.

\bigskip
Next we apply Lemma 6 and Proposition 7 to important groups  of
geometric nature. A group $G$ is called virtually  compact special
if there exists a special compact cube complex X having a finite
index subgroup of $G$ as its fundamental group (see [31] for
definition of special cube complex). The importance of virtually
compact special groups was pointed out  by Daniel Wise [31] who
proved 1-relator groups with torsion are virtually compact
special. In fact, many important groups are virtually compact
special; for example, the fundamental group of a hyperbolic
3-manifold (Agol [1]) , small cancellation groups (combination of
[31] and [1]) and hyperbolic Coxeter groups (Haglund and Wise [8])
are virtually compact
 special. Moreover,  Hangund and Wise [9] showed that quasiconvex subgroups of a
a virtually  compact special hyperbolic group $G$ (i.e.,
 a subgroups that represents  a quasiconvex subset in the  set of vertices of the Cayley graph of
$G$)  are virtual retracts of $G$. Thus the next corollary applies
in particular to this important class of subgroups.

\bigskip

\noindent {\bf 9. Corollary} {\it Let $G$ be a torsion-free
hyperbolic virtually special group and let $H$ be a quasiconvex
subgroup of $G$. Then $H$ is conjugacy distinguished.}

\bigskip
\noindent {\it Proof.} The centralizers  of torsion-free
hyperbolic groups are cyclic (cf. Proposition 12 in [22]). By
Lemma 4.1 in [21] $G$ is hereditarily  conjugacy separable. So,
using  Remark 8.2, $C_{\widehat G}(h)$ is procyclic for every
$h\in H$. Thus the result follows from Proposition
7.\hfill~\square

\bigskip
\noindent {\it Proof of Theorem B.} Combining  Theorem 1.4 in
[31], Theorem 1.2 in [11] and Proposition 4.3 in [2],  one has
that every finitely generated subgroup of $R$ is has a subgroup of
finite index which is a virtual retract of $R$. Another important
fact about $R$, proved by Newman in [23], Theorem  2 (see also
[12], p. 956), states that the centralizers of nontrivial elements
in one-relator groups with torsion are cyclic. On the other hand,
by Theorem 1.1 in [20], $R$ is hereditarily conjugacy separable.
Therefore, using Remark 8.2 above, we deduce that     ${\rm
C}_{\hat R}(g)$ is procyclic for every $g\in R$. Thus  the result
follows from Proposition 7 and the well-known fact that all
elements of finite order $n$ of $R$ are conjugate.\hfill~\square

\medskip

\noindent {\it Proof of Theorem C.} According to Theorem B in
[30], every finitely generated subgroup of $R$ is a virtual
retract of $R$. Now,  by Lemma 3.5 in [4], the centralizers of
elements of $R$ are abelian and $\overline{C_R(g)}=C_{\widehat
R}(g)$, for every $g\in R$. So the $\widehat R$-centralizers of
elements of $R$  are abelian. Thus   the result follows from
Proposition 7.\hfill~\square

\bigskip
\centerline {\bf  4. Lyndon Groups}

\bigskip

The aim of this section is to show that every finitely generated
subgroup of the Lyndon group is conjugacy distinguished. We begin
by recalling the concept of Lyndon group. We then extend the proof
in [4]  that the Lyndon group  is conjugacy separable for a more
general concept of   Lyndon group than the one considered there
(this result follows also from the result obtained by Lioutikova
[15]) . The Lyndon group was first defined in [16] with the aim of
enlarging the set of `exponents' allowed in a group. One begins
with  a free group $F$ of arbitrary rank and one wants to enlarge
$F$ to a group, usually denoted $F^{{\bf Z}[t]}$, on which the
ring of polynomials ${\bf Z}[t]$ operates (in a manner analogous
with the way the ring of integers ${\bf Z}$ operates on any
group). Myasnikov and Remeslennikov ([22])   give an explicit
construction of  the Lyndon group $F^{{\bf Z}[t]}$ as follows.

\medskip
$1$st step: One starts with the free group $F^{(0)}= F$.
\smallskip
$2$nd step: We consider  a tree of groups of the form

 $$
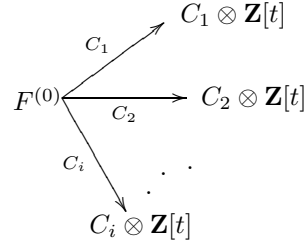
\xy\POS(0,0) *{F^{(0)} \quad}\POS(19,11) *{\quad \quad  \quad C_1\otimes {\bf Z}[t]}\POS(27,0) *{C_2\otimes {\bf Z}[t]}
\POS(18.5,-9.2) *{.}\POS(15,-10) *{.}\POS(13,-12.5)
*{.}\POS(12,-17) *{ C_i\otimes {\bf Z}[t]}
\POS(1.5,0)\ar@{->}(15,10)^{C_1}\POS(1.5,0)
\ar@{->}(18,0)_{C_2}\POS(1.5,0) \ar@{->}(10,-15)_{C_i}\endxy$$
where $\{C_i\mid 1\le i\le \delta_0\}$ is a  collection of
infinite cyclic subgroups of $F^{(0)}$  indexed by the ordinals
less than or equal a certain ordinal number $\delta_0$ and
$C_i\otimes {\bf Z}[t]=C_i\otimes_{\bf Z}{\bf Z}[t]$ is the usual
tensor product of ${\bf Z}$-modules (more precisely, the $C_i$ are
representatives of the conjugacy classes of all centralizers of
nontrivial elements of $F^{(0)}= F$, which of course in this case
are  maximal cyclic subgroups).  The edge group $C_i$ is embedded
into the vertex group $C_i\otimes {\bf Z}[t]$ by the map $C_i\to
C_i\otimes {\bf Z}[t]$ that sends $c\in C_i$ to $c\otimes 1$; this
is indeed an embedding because $C_i$ is infinite cyclic. Let
$F^{(1)}$ be the fundamental group (the tree product) of this
graph of groups. Then $F^{(1)}$ is the union of a chain
$$F=F^{(0)}= F^{(00)}\le F^{(01)}\le \cdots \le F^{(0i)}\le \cdots \le F^{(0\delta_0)}= F^{(1)},$$
where each $i$ is an ordinal, $1\le i\le \delta_0$, and if $i\ge
1$ is not a limit ordinal, then $F^{(0i)}=
F^{(0\,i-1)}*_{C_i}C_i\otimes {\bf Z}[t]$, while if $i$ is a limit
ordinal, then $F^{(0i)}= \bigcup _{j<i}F^{(0j)}$.

\smallskip
$n$th step: Here we repeat the same procedure   described in step
2, but with $F^{(0)}$ replaced with $F^{(n-1)}$ and the $C_i$
being a set of representatives of the conjugacy classes of all
cyclic centralizers of nontrivial elements of $F^{(n-1)}$ (that is
to say, intuitively,  those centralizers on which ${\bf Z}[t]$
does not operate yet).
\smallskip

Then the Lyndon group is
$$F^{{\bf Z}[t]}= \bigcup_{m=0}^\infty F^{(m)}$$
 (cf.[22], Theorem 8).  Therefore one can describe $F^{{\bf Z}[t]}$ as a union of a chain of groups
$$F= F_{(0)}\le  F_{(1)}\le \cdots \le  F_{(i)}\le \cdots \le  F_{(\delta)}= F^{{\bf Z}[t]}$$
indexed by the ordinals $i$ less than or equal to   a certain
ordinal $\delta$ and such that $F_{(i)}=
F_{(i-1)}*_{C_i}(C_i\otimes {\bf Z}[t])$, when $i$ is a nonlimit
ordinal, while if $i$ is  limit ordinal, then $F_{(i)}= \bigcup
_{j<i}F_{(j)}$.

Observe that the image $A_i=C_i\otimes 1$ of $C_i$ in $C_i\otimes
{\bf Z}[t]$ is a direct summand of $C_i\otimes {\bf Z}[t]$. Say
$C_i\otimes {\bf Z}[t]= A_i\oplus \tilde A_i$, where $A_i\cong
C_i\cong {\bf Z}$ and $\tilde A_i$ is a free abelian group of
infinite rank.  We identify $A_i$ with $C_i$.

\bigskip
\noindent {\bf 10. Lemma } {\it Let $0\le i\le s \le \delta$. Then
\medskip

\item{\rm (a)} there exists a group epimorphism $\varphi_{s,i}:
F_{(s)}\rightarrow F_{(i)}$ which is the identity on the subgroup
$F_{(i)}$  of $F_{(s)}$;
\smallskip
\item{\rm (b)} $F_{(s)}= K_{s,i}\semidirprod F_{(i)}$, for some
normal subgroup $K_{s,i}$ of $F_{(s)}$. }
\medskip
 \noindent {\bf  Proof.}
\medskip
(a) We shall use transfinite induction to define a group
homomorphism $\varphi_{s,i}$ from $F_{(s)}$ onto $F_{(i)}$ which
is the identity on the subgroup $F_{(i)}$ of $F_{(s)}$, and such
that the restriction of $\varphi_{s,i}$ to $F_{(k)}$ is
$\varphi_{k,i}$  whenever $i\le k\le s$. Define $\varphi_{i,i}$ to
be the identity map and assume that $\varphi_{r,i}$ has already
been  defined for all $r<s$ ($i\le r\le s\le \delta$). Then we
define $\varphi_{s,i}: F_{(s )}\rightarrow F_{(i)}$ as follows:
\medskip
-- if $s$ is a limit ordinal, put $\varphi_{s,i}= \bigcup_{r<
s}\varphi_{r,i}$, and

\smallskip

-- if $s$ is a nonlimit ordinal, then define first
$$\psi: F_{(s )}= F_{(s -1)}*_{C_{s -1}} (C_{s -1}\oplus \tilde A_{s -1})\rightarrow F_{(s -1)}$$
 by sending $F_{(s -1)}$ identically to  $F_{(s -1)}$ and sending $\tilde A_{s -1}$ to $1$. Then define  $\varphi_{s,i}=  \varphi_{s -1, i} \psi$.
\smallskip
(b) This is clear from (a).\hfill~\square

 \bigskip

 The next task is to show that $F^{{\bf Z}[t]}$ is a conjugacy separable group. In fact we will prove more generally that  the Lyndon group belongs to a class ${\cal X}$
 of abstract groups  that  satisfy a series of properties including that of being conjugacy separable. This class ${\cal X}$ was introduced in [28] and we describe it briefly here.
An  abstract group $R$ is in ${\cal X}$   if
\medskip
\item{\rm (a)} $R$ is conjugacy separable (so that  in particular
$R\le \hat R$);
\smallskip
\item{\rm (b)} $R$ is quasi-potent (i.e., for every cyclic
subgroup $H$ of $R$, there exists a subgroup $K$ of finite index
in $H$ such that every subgroup of finite index of $K$ is of the
form $K\cap N$, for some normal subgroup $N$ of finite index in
$R$);
\smallskip
\item{\rm (c)} whenever $A$ and $B$  are cyclic subgroups  of $R$,
the set $AB$ is closed in the profinite topology of $R$;
\smallskip
\item{\rm (d)} every cyclic subgroup of $R$ is conjugacy
distinguished, i.e., if $C$ is a cyclic subgroup of $R$  and
$a\in R$, then $a^R\cap C= \emptyset$ if and only if $a^{\hat
R}\cap \bar C = \emptyset$;
\smallskip
\item{\rm (e)}  if $A$ and $B$ are cyclic subgroups of $R$, then
$A\cap B=1$ if and only if $\bar A \cap \bar B=1$; and
\smallskip
\item{\rm (f)}  if $A= \langle a \rangle$ is an infinite cyclic
subgroup of $R$, and $\gamma\in \hat R$ with  $\gamma \in {\rm
N}_{\hat R}(\bar A)$, then $\gamma \in {\rm N}_{\hat R}(A)$, i.e.,
$\gamma a \gamma^{-1}\in \{a, a^{-1}\}$.

 \bigskip

\noindent {\bf 11. Proposition}   {\it The Lyndon group $F^{{\bf
Z}[t]}$ is in the class $\cal X$, and in particular it is
conjugacy separable.}

\medskip
 \noindent {\bf  Proof.}     We continue with the above notation. We   shall prove inductively that in  fact  each $F_{(s )}$  is in the class $\cal X$,    for all $0\le s \le \delta$.
It is well-known that the free group  $F_{(0)}= F$  is in class
$\cal X$. Assume that  $F_{(j)}\in {\cal X}$  for $0\le j<s $. If
$s$ is a nonlimit ordinal, then $F_{(s )}= F_{(s -1)}*_{C_{s -1}}
(C_{s -1} \otimes {\bf Z}[t])$ is in $\cal X$   according to
Theorem A  in [28], since $C_{s -1} \otimes {\bf Z}[t]$ is free
abelian of infinite rank, and so both $C_{s -1} \otimes {\bf
Z}[t]$ and $F_{(s -1)}$ are in the class $\cal X$.

Let now  $s$ be a limit ordinal.  Then $F_{(s )}=
\bigcup_{j<s}F_{(j)}$. We have to verify that $F_{(s )}$ satisfies
properties (a)-(f) of class $\cal X$. Observe first that $F_{(s
)}$ is residually finite: indeed, let $1\not= x\in F_{(s )}$; then
$x\in F_{(i)}$, for some $i<s $, and since $F_{(i)}$ is residually
finite, there exists some   $N\triangleleft_f F_{(i)}$ with
$x\not\in N$; so, if $\varphi_{s,i}: F_{(s)}\rightarrow F_{(i)}$
denotes the epimorphism defined in   Lemma 10, we have that
$\varphi_{s,i}^{-1}(N)$ is a normal subgroup of finite index in
$F_{(s)}$ that misses $x$. It follows from this and part (b) of
Lemma 10 that $F_{(i)}$ is closed in the profinite topology of
$F_{(s)}$, and moreover  the profinite topology of $F_{(s)}$
induces on $F_{(i)}$ its full profinite topology (cf. Lemma 3.1.5
in [26]). Using this one easily verifies that   $F_{(s)}$
satisfies  properties (c)  and (e) of class $\cal X$.

To verify property (a) of class $\cal X$ (i.e., that  $F_{(s)}$ is
conjugacy separable) let $x,y\in F_{(s)}$ and assume that  $y=
x^\gamma= \gamma^{-1}x\gamma$, where $\gamma\in \widehat
{F_{(s)}}$. Then there exists some ordinal $i$, $0\le i< s $ with
$x,y\in F_{(i)}$.  Let $\varphi_{s,i}: F_{(s)}\rightarrow F_{(i)}$
be an epimorphism such that $\varphi_{s}$ is the identity on
$F_{(i)}$ (see Lemma 10). Let $\widehat {\varphi_{s,i}}: \widehat
{F_{(s)}}\rightarrow \widehat {F_{(i)}}$ be the continuous
homomorphism induced by $\varphi_{s,i}$ (cf. [26], Lemma 3.2.3).
Put $\tilde \gamma= \widehat {\varphi_{s,i}} (\gamma)$. Then $y=
x^{\tilde \gamma}$. Since, by assumption,  $F_{(i)}$ is conjugacy
separable, there exists some $c\in F_{(i)}\le F_{(s)}$ such that
$y=x^c$, as needed.

To verify property (b) of class $\cal X$ (i.e., that  $F_{(s)}$ is
quasi-potent) let  $H$ be a cyclic subgroup of  $F_{(s)}$.  Then
$H\le  F_{(i)}$ for some ordinal $i<s$. By Lemma 10, $F_{(s)}=
K_{s,i}\semidirprod F_{(i)}$. Since $F_{(i)}$ is quasi-potent,
there exists a subgroup $K$    of  $H$ of finite index such that
every subgroup of finite index of $K$  has the form $K\cap U$, for
some  $U\triangleleft_f F_{(i)}$. Since $K\cap K_{s,i}U= K\cap U$
and $K_{s,i}U\triangleleft _f F_{(s)}$, we deduce that $F_{(s)}$
is quasi-potent.

For property (d) of class $\cal X$, let $C$ be a cyclic subgroup
of $F_{(s)}$ and let $a^{F_{(s)}}\cap C=\emptyset$, where $a\in
{F_{(s)}}$. Assume that $a^\gamma= \alpha\in \bar C$, for some
$\gamma\in \widehat {{F_{(s)}}}$, where $\bar C$ is the closure of
$C$ in   $ \widehat {{F_{(s)}}}$. Let $i$ be an ordinal, $0\le
i<s$ such that $C\le F_{(i)}$, let $\varphi_{s,i}:
F_{(s)}\rightarrow F_{(i)}$ be the epimorphism described in Lemma
10,  and  let $\hat \varphi_{s,i}: \widehat {F_{(s)}}\rightarrow
\widehat {F_{(i)}} = \overline {F_{(i)}}$ be the induced
epimorphism. Note that $\bar C$ is also the closure of $C$ in
$\overline {F_{(i)}}$, and hence $\hat \varphi_{s,i} (\bar C)=
\bar C$. Put $\tilde \gamma= \hat \varphi_{s,i}(\gamma)$ and
$\tilde \alpha= \hat \varphi_{s,i}(\alpha)$. Then $a^{\tilde
\gamma}= \tilde \alpha\in \bar C$. Since  $F_{(i)}$ has property
(d), there exists $c\in F_{(i)}\le F_{(s)}$ with $a^c\in C$, a
contradiction. Hence  $a^{\overline {F_{(s)}}}\cap \bar
C=\emptyset$, showing that
 $F_{(s)}$ has property (d).

Finally, we check that $F_{(s)}$ has property (f) of class $\cal
X$. Let $A= \langle a\rangle$ be an infinite cyclic subgroup of
$F_{(s)}$ and assume that $\gamma^{-1}a \gamma \in  \bar A$, for
some $\gamma \in \widehat {F_{(s)}}$. Let $i$ be an ordinal, $0\le
i<s$ such that $a \in F_{(i)}$, let $\varphi_{s,i}:
F_{(s)}\rightarrow F_{(i)}$ be the epimorphism described in Lemma
10 and let $\hat \varphi_{s,i}: \widehat {F_{(s)}}\rightarrow
\widehat {F_{(i)}} = \overline {F_{(i)}}$ be the induced
epimorphism.  Put $\tilde \gamma= \hat \varphi_{s,i}(\gamma)$ and
observe that $\tilde \gamma^{-1} a \tilde \gamma \in  \bar A$ in
$\overline  {F_{(i)}}$, since $\bar A$ is the closure of $A$ in
both  $F_{(s)}$  and $F_{(i)}$. Since $F_{(s)}$ has property (f),
we have that $\tilde \gamma^{-1} a \tilde \gamma$ is either $a$ or
$a^{-1}$.  Since $\hat \varphi_{s,i}$ is the identity on $\bar A$,
we deduce that also $\gamma^{-1} a\gamma\in \{a,
a^{-1}\}$.\hfill~\square

\bigskip

 \noindent {\bf 12. Proposition} {\it Let $H= \langle h_1, \dots, h_n\rangle$ be  a finitely generated subgroup of
the Lyndon group $F^{{\bf Z}[t]}$. Then
\medskip
\item{\rm (a)}  there exists a finitely generated subgroup $K$  of
$F^{{\bf Z}[t]}$ such that $H\le K$  and $K$ is a retract of
$F^{{\bf Z}[t]}$; and
\smallskip
\item{\rm (b)}
 $H$ is a virtual retract of  $ F^{{\bf Z}[t]}$. }

\medskip

\noindent {\it Proof.}  We continue with the above description and
notation. Since $H$ is finitely generated, there exists a smallest
ordinal $s$ with $s<\delta$ such that $H\le F_{(s)}$. Note that
$s$ cannot be a limit ordinal.
\medskip
(a)  We prove this by induction on $s$. If $s=0$, this
 is clear since $F_{(0)}$ is a free group.
Suppose $s>0$ and that for every ordinal $i$, with $i<s$, every
finitely generated subgroup of $F_{(i)}$ is contained in a
finitely generated retract   of $F_{(i)}$. Since $s$ is not a
limit ordinal, one has
$$F_{(s)}= F_{(s-1)}*_C(C\otimes {\bf Z}[t]).$$
Write each $h_i$ as a product of elements of    $F_{(s-1)}$ and
$C\otimes {\bf Z}[t]$. Say  $f_1, \ldots, f_t$ are all the
elements of $F_{(s-1)}$ involved in those products, and let $b_1,
\dots, b_r$ be all the elements of $C\otimes {\bf Z}[t]$ involved
in those products. By hypothesis   there exists a finitely
generated retract $K_1$ of $F_{(s-1)}$ such that  $K_1\ge \langle
f_1, \dots, f_t, C \rangle$. Let $B$ be a finitely generated
direct summand of the free abelian group $C\otimes {\bf Z}[t]$
such that $b_1, \dots, b_r\in B$. Then  $H\le K_1*_CB$,  $K_1*_CB$
is finitely generated and  a retract of $F_{(s)}$. This proves
part (a).
\medskip
(b) Let $K$ be as in part (a). Since $K$ is a limit group and
$H\le K$, there is a subgroup $U$ of finite index in $K$ with
$U\ge H$ and an epimorphism $\varphi: U\to H$ which is the
identity on $H$   (cf. [30], Theorem B).  Let $\psi: F^{{\bf
Z}[t]} \to K$ be a retraction and put $V= \psi^{-1}(U)$. Then $V$
has finite index in $F^{{\bf Z}[t]}$ and  the composite $\varphi
\psi_{|V}:V\to H$ is a retraction, as needed.\hfill~\square

\bigskip
\noindent {\it Proof of Theorem D.}
Let $H$ be a finitely generated subgroup of
$F^{{\bf Z}[t]}$ and assume that $g^\gamma\in \bar H$, where    $g\in F^{{\bf Z}[t]}$ and
 $\gamma\in \widehat{F^{{\bf Z}[t]}}$. One needs to show that there exists some $c\in F^{{\bf Z}[t]}$ with $g^c\in H$.
 By Proposition 12 there
exists a retraction  $\varphi: F^{{\bf Z}[t]}\to K$, where $K$  is a finitely generated subgroup of  $F^{{\bf Z}[t]}$
containing $\langle H, g\rangle$.  Then $\varphi$ extends to a retraction $\hat \varphi: \widehat{F^{{\bf Z}[t]}}\to \hat K= \bar K$.
Put $\gamma'= \hat \varphi (\gamma)$. Then $g^{\gamma'}\in \bar H$.
Since every finitely generated
subgroup of $F^{{\bf Z}[t]}$ is a limit group,   the result follows from
Theorem C. \hfill~\square

\vskip 1cm

\centerline {\bf REFERENCES}
\bigskip

\noindent [1] I. Agol,  The virtual Haken conjecture. {\it Doc.
Math.}, {\bf 18} (2013) 1045--1087.  With an appendix by Agol,
Daniel Groves, and Jason Manning.

\smallskip
\noindent [2] M. Bestvina, Geometric group theory and 3-manifolds
hand in hand: the fulfillment of Thurston's, {\it Bull. Amer.
Math. Soc.} {\bf 51},  53--70 (2014).

\smallskip
\noindent [3] O. Bogopolski and K.-U. Bux, Subgroup conjugacy
separability for surface groups, arXiv:1401.6203.

\smallskip
\noindent [4]  S. C. Chagas and  P. A. Zalesskii, Limit groups are
conjugacy separable, {\it Int. J. Algebra Comput.} {\bf 17},
851--857 (2007).

\smallskip
\noindent  [5]  J. Dyer, Separating conjugates in amalgamated free
products and HNN extensions, {\it J. Austral. Math. Soc.} (Series
A) {\bf 29}, 35--51 (1980).
\smallskip

\noindent  [6] J.L. Dyer and P. Scott,  Periodic automorphisms of
free groups, {\it Comm. Alg.} {\bf 3},  195--201 (1975).

\smallskip

\noindent [7] M. Ferov, On conjugacy separability of graph
products of groups.  arXiv:1409.8594.

\smallskip

\noindent [8] F. Haglund, D.T. Wise, Coxeter groups are virtually
special. {\it Advances in Mathematics} {\bf 224} (2010)
1890--1903.

\smallskip
 \noindent [9]  F. Haglund, D.T. Wise, Special cube
complexes. {\it Geom. Funct. Anal.} {\bf 17} (2008), no. 5,
1551--1620.

\smallskip
\noindent [10] M. Hall Jr., Coset representations in free groups,
{\it Trans. Amer. Math. Soc.} {\bf 67}, 421--432 (1949).

\smallskip

\noindent [11] G. C. Hruska and D. T. Wise, Towers, ladders and
the B. B. Newman spelling theorem, {\it J. Aust. Math. Soc.} {\bf
71}, 53--69 (2001).

\smallskip

\noindent [12] A. Karrass, A.  and D. Solitar,  The free product
of two groups with a malnormal amalgamated subgroup. {\it Canad.
J. Math.} {\bf 23},  933--466 (1971).

\smallskip

\noindent [13] A. Karrass, A. Pietrowski and D. Solitar, Finite
and infinite cyclic extensions of free groups, {\it J. Austral.
Math. Soc.} {\bf 16}, 458--466 (1973).

\smallskip

\noindent [14] O. Kharlampovich and A. Myasnikov, Irreducible
affine varieties over a free group. II. Systems in triangular
quasi-quadratic form and description of residually free groups.
{\it  J. Algebra} {\bf 200}, 517--570 (1998).

\smallskip
\noindent [15]  E. Lioutikova, Lyndon's group is conjugately
residually free, {\it Int.  J. Algebra Comput.} {\bf 13}, 255--275
(2003).

\smallskip

\noindent [16] R. Lyndon, Groups with parametric exponents, {\it
Trans. AMS} {\bf 96}, 518-533 (1960).

\smallskip
\noindent [17] 8. R.C.  Lyndon  and P.E. Schupp,   {\it
Combinatorial Group Theory}, Springer,  Heidelberg, 1977.
\smallskip

\noindent [18] W.  Magnus,  A. Karrass and D. Solitar, {\it
Combinatorial Group Theory}, J. Wiley and Sons, New York, 1966.

\smallskip
\noindent [19] A. Minasyan, Hereditary conjugacy separability of
right angled Artin groups and its applications, {\it Groups Geom.
Dyn.} {\bf 6}, 335-388 (2012).


\smallskip
\noindent [20] A. Minasyan and P. Zalesskii, One-relator groups
with torsion are conjugacy separable, {\it J. Algebra} {\bf 382},
39-45 (2013).

\smallskip
\noindent [21] A. Minasyan and P. Zalesskii, Virtually compact
special hyperbolic groups are conjugacy separable.
arXiv:1504.00613.

\smallskip

\noindent [22]  A.G. Myasnikov and V.N. Remeslennikov, Exponential
groups 2 extensions of centralizers and tensor completion of
CSA-groups,  {\it Int.  J. Algebra and Comput.}, {\bf 6} no. 6,
687-711 (1996).

\smallskip
\noindent [23] B.B. Newman, Some results on one-relator groups.
{\it Bull. Amer. Math. Soc.} {\bf 74}, 568--571 (1968).

\smallskip
\noindent [24] L. Ribes and P. Zalesskii, The pro-$p$ topology of
a free group and algorithmic problems in semigroups, {\it Int. J.
Algebra and Computation} {\bf 4}, 359--374 (1994).

\smallskip
\noindent [25] L. Ribes and P. Zalesskii, Conjugacy separability
of amalgamated free products of groups, {\it J. Algebra} {\bf
179}, 751--774 (1996).

\smallskip

\noindent [26] L. Ribes and P. Zalesskii, {\it Profinite Groups},
2nd edition,  Springer, Berlin, 2010.

\smallskip

\noindent [27] L. Ribes and P. Zalesskii, Normalizers in groups
and in their profinite completions, {\it Rev. Mat. Iberoam.} {\bf
30}, 165--190 (2014).

\noindent [28] L. Ribes, D. Segal and P.A. Zalesskii,
 Conjugacy separability and free products of groups with cyclic amalgamation, {\it J. London Math. Soc.} {57}, 609--628 (1998).

\smallskip
\noindent [29] G.P. Scott,   An embedding theorem for groups with
a free subgroup of finite index. {\it Bull. London Math. Soc.}
{\bf 6}, 304--306 (1974).

\smallskip
\noindent [30] H. Wilton, Hall's Theorem for limit groups, {\it
GAFA} {\bf 18}, 271-303 (2008).

\smallskip

\noindent [31] D.T. Wise, The structure of groups with a
quasiconvex hierarchy, {\it Electr. Res. Announcement
 in Math. Sci.}, Volume 16,
 44--55 (2009).

\smallskip

\noindent  [32] P.A. Zalesskii and O.V. Mel'nikov, Subgroups of
profinite groups acting on trees, {\it Math. USSR Sbornik}, {\bf
63}, 405--424 (1989).

\vskip1cm
\noindent  Luis Ribes: School of Mathematics and Statistics, Carleton University, Ottawa ON K1S 5B6, Canada

\noindent  E-mail: lribes@math.carleton.ca

\medskip
\noindent  Pavel Zalesskii: Departamento de Matem\'atica, Universidade de Bras\'{\i}lia, 70910-900
Bras\'{\i}lia-DF, Brazil.

\noindent  E-mail: pz@mat.unb.br
  \end